\documentclass[a4paper,11pt]{article}
\usepackage{amssymb}
\usepackage{amssymb, color, eufrak}
\usepackage{amsmath}
\usepackage{amsthm}

\newcommand{\aand}{\quad\mbox{and}\quad}

\newcommand{\R}{\mathbb{R}}

\newcommand{\beq}{\begin{equation} }
\newcommand{\eqq}{\end{equation} }
\newcommand{\cuad}{{\sqcap\kern-.68em\sqcup}}

\newcommand{\norm}[1]{\|#1\|}

\newtheorem{teo}{Theorem}[section]

\newtheorem{proposition}{Proposition}[section]

\newtheorem{lemma}{Lemma}[section]

\newtheorem{remark}{Remark}[section]
\newcommand{\bremark}{\begin{remark} \em}
\newcommand{\eremark}{\end{remark} }

\def \l{\lambda}
\def \e{\varepsilon}
\def \ds{\displaystyle}

\def\beeq{\begin{equation}}
\def\eeq{\end{equation}}
\newcommand{\begeqaet}{\begin{align*}}
\newcommand{\eneqaet}{\end{align*}}

\begin{document}
\begin{center}{\bf  \Large
On Gaussian curvature equation in $\R^2$  with \\[2mm]

prescribed nonpositive curvature
  }
%%%%%%%%%%%%%%%%%%%%%%%%%%%%%%%%%%%%%%%%%%%%%%%%%%%%%%%%%%%%%%%%%%%%%%
%%%%%%%%%%%%%%%%%%%%%%%%%%%%%%%%%%%%%%%%%%%%%%%%%%%%%%%%%%%%%%%%%%%%%%
\bigskip\bigskip
{\small

{\sc  Huyuan Chen\footnote{chenhuyuan@yeah.net}}

\medskip
Department of Mathematics, Jiangxi Normal University,\\
Nanchang, Jiangxi 330022, PR China \\[10pt]
 {\sc  Dong Ye\footnote{dye@math.ecnu.edu.cn, dong.ye@univ-lorraine.fr}}

\medskip
Center for PDEs, School of Mathematical Sciences, East China Normal University,\\
Shanghai Key Laboratory of PMMP, Shanghai 200062, PR China

\smallskip
IECL, UMR 7502, University of Lorraine, 57050 Metz, France \\[12pt]
{\sc  Feng Zhou\footnote{fzhou@math.ecnu.edu.cn}}

\medskip
Center for PDEs, School of Mathematical Sciences, East China Normal University,\\
Shanghai Key Laboratory of PMMP, Shanghai 200062, PR China \\[20pt]
}
\end{center}

\begin{center}
{\sl Dedicated to Professor Weiming Ni's seventieth birthday}
\end{center}

\begin{abstract}
The purpose of this paper is to study the solutions of
$$
\Delta u +K(x) e^{2u}=0 \quad{\rm in}\;\; \mathbb{R}^2
$$
with $K\le 0$. We introduce the following quantity:
$$\alpha_p(K)=\sup\left\{\alpha \in \R:\, \int_{\R^2} |K(x)|^p(1+|x|)^{2\alpha p+2(p-1)} dx<+\infty\right\}, \quad \forall\; p \ge 1.$$
Under the assumption $({\mathbb H}_1)$: $\alpha_p(K)> -\infty$ for some $p>1$ and $\alpha_1(K) > 0$, we show that for any $0 < \alpha < \alpha_1(K)$, there is a unique solution $u_\alpha$ with $u_\alpha(x) = \alpha \ln |x|+ c_\alpha+o\big(|x|^{-\frac{2\beta}{1+2\beta}} \big)$ at infinity and $\beta\in (0,\,\alpha_1(K)-\alpha)$.
Furthermore, we show an example $K_0 \leq 0$ such that $\alpha_p(K_0) = -\infty$ for any $p>1$ and $\alpha_1(K_0) > 0$, for which we study the asymptotic behavior of solutions. In particular, we prove the existence of a solution $u_*$ such that $u_* -\alpha_*\ln|x| = O(1)$ at infinity for some $\alpha_* > 0$, but who does not converge to a constant at infinity. This example exhibits a new phenomenon of solutions with logarithmic growth and non-uniform behavior at infinity.
\end{abstract}



%%%%%%%%%%%%%%%%%%%%%%%%%%%%%%%%%%%%%%%%%%%%%%%%%%%%%%%%%%%%%%%%%%%%%%
%%%%%%%%%%%%%%%%%%%%%%%%%%%%%%%%%%%%%%%%%%%%%%%%%%%%%%%%%%%%%%%%%%%%%%
%\keywords{}
%\subjclass{}

\setcounter{equation}{0}
\section{Introduction}
In this paper, we study the solutions of
\begin{equation}\label{eq 1.1}
\Delta u + K(x)e^{2u}=0 \quad{\rm in}\;\; \mathbb{R}^2,
\end{equation}
where $K$ is a nonpositive, non trivial and locally bounded function in $\R^2$.

\medskip
Equation (\ref{eq 1.1}) finds its origin from conformal geometry:  Given a Riemannian surface $(\mathcal{M}, g)$ and a function $K$ on $\mathcal{M}$, we may ask whether there is a new Riemannian metric $g_1$, conformal to $g$, such that $K$ is the Gaussian curvature of $g_1$. In other words, we need to find a function $u$ on $\mathcal{M}$ such that $g_1=e^{2u} g$  and $u$ is a solution of the elliptic equation
\begin{equation}\label{eq 1.2}
  \Delta_g u - k_g +Ke^{2u}=0\quad{\rm on}\;\; \mathcal{M}.
\end{equation}
Here $\Delta_g$ is the Laplace-Beltrami operator, $k_g$ is the Gaussian curvature on $({\mathcal M}, g)$. In the Euclidean case $\mathcal{M}=\R^2$ and $g = |dx|^2$,
(\ref{eq 1.2}) reduces to equation (\ref{eq 1.1}).

\medskip
The equation (\ref{eq 1.1}) has been studied extensively, especially for the nonpositive, nonconstant curvature case, see \cite{K, S, N, McO1, CN1} and the references therein.
The nonexistence of entire solution in $\R^2$ with
$K \equiv -1$ was given by Ahlfors \cite{Ah} in 1938, and then improved by Sattinger \cite{S}, Ni \cite{N} for
$K \leq 0$ in $\R^2$ and $K(x)\le - C|x|^{-2}$ at infinity with $C > 0$. The first general existence result for (\ref{eq 1.1}) is due to Ni \cite{N}, he proved that if $K\le 0$, $K\not\equiv 0$ in $\R^2$ with $|K(x)| \le C|x|^{-\ell}$ at infinity for some $\ell > 2$, then for any $\alpha \in (0, \beta)$
with $\beta = \min \{8, \frac{\ell-2}{3}\}$, there exists a solution $u$ satisfying $u(x)= \alpha \ln |x| + O(1)$ at infinity. Later on, it was proved by Cheng-Ni \cite{CN1} that
if $K \leq 0$ in $\R^2$, $ K(x) \sim  - |x|^{-\ell}$ at infinity with $\ell>2$, then for any $\alpha\in (0, \,  \frac{\ell-2}{2})$, (\ref{eq 1.1}) has a unique solution satisfying
\begin{align}
\label{newb1}
u_\alpha(x) = \alpha\ln|x| + c_\alpha + o(1) \;\; \mbox{as }\; |x|\to \infty, \quad \mbox{with } c_\alpha \in \R.
\end{align}
The solutions set has layer structure, that is $u_\beta < u_\alpha$ in $\R^2$, for any $0 < \beta < \alpha <  \alpha_1(K)$.
Moreover, a unique maximal solution $U$ exists such that $U(x) = \frac{\ell-2}{2}\ln |x|-\ln\ln |x|+O(1)$ at infinity.

\medskip
For general curvature functions $K$, the following quantity was introduced by Cheng-Ni in \cite{CN2},
\begin{equation}\label{alpha 1-0}
\alpha_1(K)=\sup\left\{\alpha \in \R:\, \int_{\R^2} |K(x)|(1+|x|)^{2\alpha} dx< +\infty\right\}
\end{equation}
with the convention $\sup\emptyset = -\infty$. When $K$ is controlled polynomially at infinity, using $\alpha_1(K)$,
Cheng-Ni characterized the existence to problem \eqref{eq 1.1}: Assume that $K \leq 0$ satisfies
$$\mbox{$\alpha_1(K) > 0$ and $|K(x)|< C(|x|^m + 1)$ in $\R^2$ for some $m>0$.} \leqno{( \mathbb{H}_0)}$$
Then for any $\alpha\in (0, \alpha_1(K))$, \eqref{eq 1.1} has a unique solution such that
 \begin{equation}\label{b1}
 u_\alpha = \alpha\ln |x| + O(1) \quad{\rm as}\ \ |x|\to\infty.
 \end{equation}

Later on, by employing a variational method with weighted Sobolev spaces (firstly applied by McOwen in \cite{McO}),
Cheng-Lin \cite{CL} succeeded to handle the existence of problem \eqref{eq 1.1} without any growth control for $K$:
Assume that $\alpha_1(K) > 0$ for a nonpositive, locally bounded function $K$.
\begin{itemize}
\item For any $\alpha \in (0, \alpha_1(K))$, equation \eqref{eq 1.1} possesses a unique solution $u_\alpha$ satisfying
\begin{align}
\label{ChenL}
\int_{\R^2}  Ke^{2u}  dx = - 2\alpha \pi \aand \int_{\R^2\backslash B_1} \Big|\nabla(u - \alpha\ln|x|)\Big|^2 dx<+ \infty.
\end{align}
\item Moreover, for any $0 < \beta < \alpha <  \alpha_1(K)$, there holds $u_\beta < u_\alpha$ in $\R^2$.
\item If $\alpha_1(K) < \infty$, then $\lim_{\alpha \to \alpha_1(K)} u_\alpha = u_{\alpha_1}$ is a solution to \eqref{eq 1.1} satisfying \eqref{ChenL} with $\alpha = \alpha_1(K)$.
\end{itemize}

The quantity
\begin{equation*}
 \mathcal{C}(u) := \int_{\R^2} K  e^{2u} dx
 \end{equation*}
 is called the total curvature of the conformal metric $g_u = e^{2u}|dx|^2$. It plays an important role in the analysis of \eqref{eq 1.1} and has a notable meaning in geometry. There are many interesting discussions in \cite{CL}. For example, by Theorem 1.4 in \cite{CL}, when $K \leq 0$, a solution $u$ to \eqref{eq 1.1} is logarithmically up bounded at infinity if and only if ${\mathcal C}(u)$ is finite and the metric $e^{2u}|dx|^2$ is complete.

\medskip
However, to claim the asymptotic behavior like \eqref{b1}, Cheng-Lin need additional assumptions like $({\mathbb H}_0)$. A naturel question is to understand the asymptotic behavior of solutions provided by variational method in \cite{CL}. We will point out that any solution satisfying \eqref{b1} is just a Cheng-Lin type solution, see Proposition \ref{lnew1} below. Moreover, we want to prove the existence of solutions satisfying \eqref{b1} without any growth assumption on the curvature function $K$.

\medskip
Another naturel question is to ask if there exist solutions with logarithmic growth but anisotropic behavior at infinity, more precisely, solutions of \eqref{eq 1.1} satisfying \eqref{b1} but not \eqref{newb1}. Our purpose here is to answer these questions.

\medskip
In fact,  we propose to replace the pointwise growth assumption $(\mathbb{H}_0)$ by an integral estimate. To this end, consider the following quantity related to $K$, which generalizes $\alpha_1(K)$:
 \begin{equation}\label{alpha 1}
 \alpha_p(K)= \sup\left\{\alpha \in \R:\, \int_{\R^2} |K(x)|^p(1+|x|)^{2\alpha p+2(p-1)} dx< +\infty\right\}, \quad \forall\; p \ge 1.
 \end{equation}
We use always the convention $\sup\emptyset = -\infty$. Our first main result states as follows.
\begin{teo}\label{teo 1}
Let $K\le 0$ be locally bounded in $\R^2$. Assume that
$$
 \mbox{$\alpha_1(K) > 0$ and there exists $p>1$ such that $\alpha_p(K) > -\infty$.} \leqno{(\mathbb{H}_1)}
$$
Then for any $\alpha\in(0,\,\alpha_1(K))$, problem (\ref{eq 1.1}) possesses a unique solution $u_\alpha$
satisfying \eqref{b1}. Moreover, the total curvature of $u_\alpha$ is $-2\alpha\pi$; and there exists $c_\alpha \in \R$ such that for any $\beta\in (0, \alpha_1(K)-\alpha)$,
\begin{equation}\label{1.3}
 u_\alpha(x)=\alpha \ln|x| +c_\alpha+o\left(|x|^{-\frac{2\beta}{1+2\beta}}\right)\quad {\rm as }\; |x| \to \infty.
\end{equation}
\end{teo}

\begin{remark}
Combining with Proposition \ref{lnew1} below, Theorem \ref{teo 1} means that under $(\mathbb{H}_1)$, for any $\alpha\in(0,\,\alpha_1(K))$, the unique solution satisfying \eqref{ChenL} given by \cite{CL} has the asymptotic behavior \eqref{b1}.
\end{remark}

\begin{remark}\label{re 1.1}
Assumption $(\mathbb{H}_1)$ is weaker than $(\mathbb{H}_0)$, see Lemma \ref{pr 2.3} below.
\end{remark}

Our construction is based on the understanding of $\alpha_p(K)$. Indeed, under $({\mathbb H}_1)$, we can prove that $\alpha_p(K)$ are positive for $p > 1$ but nearby. This observation allows us to construct suitable super and sub solutions of (\ref{eq 1.1}). Let $w_0$ be a positive smooth radial function such that
\begin{equation}\label{defw}
 w_0(x) = \ln |x|\quad\ {\rm in}\ \, \R^2\setminus B_1,
\end{equation}
where $B_1 = \{x \in \R^2, |x| < 1\}.$
In the following, $B_r(x)$ stands for the standard open ball with center $x\in \R^2$ and the radius $r$; $B_r$ denotes always $B_r(0)$. We will look for a solution $u_\alpha$ to (\ref{eq 1.1})  in the form
$$
 u_\alpha =  \alpha  w_0 + v_\alpha,
$$
where $v_\alpha$   resolves
\begin{equation*}
- \Delta v  = Ke^{2\alpha w_0} e^{2v} + \alpha \Delta w_0 \quad{\rm in}\;\; \mathbb{R}^2.
\end{equation*}
The crucial point is to know whether the above equation has a bounded solution, or equivalently, whether $Ke^{2\alpha w_0} e^{2v} + \alpha \Delta w_0$ belongs to the Kato's class. If the answer is affirmative, the Perron's method can be employed to get a solution of (\ref{eq 1.1}).

\medskip
Our second concern is to understand solutions to \eqref{eq 1.1} with an exotic $K_0 \leq 0$ such that
\begin{equation}
\label{aK0}
\alpha_p(K_0)=-\infty\;\; \mbox{for any }\; p> 1 \aand \alpha_1(K_0) > 0.
\end{equation}
We construct here an explicit model $K_0$ as follows. Consider the sequence $a_n := (n, 0)_{n\ge 2}\in \R^2$. Let
\begin{equation}\label{k2}
 K_0(x)=-\sum^\infty_{n=2} r_n^{-2}n^{-\ell}\eta_0\left(\frac{x-a_n}{r_n}\right) \quad \mbox{in } \R^2,
\end{equation}
where
\begin{align*}
\ell > q > 1, \quad r_n = e^{-n^q} \;\; \mbox{for all }\;  n \geq 2
\end{align*}
and $\eta_0$ is a smooth radial cut-off function satisfying
 \begin{equation}\label{eta 0}
\chi_{B_\frac12} \leq \eta_0 \leq \chi_{B_1}.
\end{equation}
Here $\chi_\Omega$ means the usual characteristic function of $\Omega$.
Observe that $r_n \leq e^{-2} < \frac{1}{4}$ for any $n \geq 2$, so
\begin{align*}
\int_{\R^2} |K_0(x)|^p(1+|x|)^{2\alpha p+2(p-1)} dx = \sum^\infty_{n=2} \int_{B_{r_n}(a_n)} |K_0(x)|^p(1+|x|)^{2\alpha p+2(p-1)} dx =: \sum^\infty_{n=2}  I_n.
\end{align*}
Direct computation implies then: for $p\ge 1$, $\alpha \in \R$, (as $n$ goes to $\infty$)
\begin{align*}
I_n & := r_n^{-2p} n^{\ell p}\int_{B_{r_n}(a_n)}(1+ |x|)^{2\alpha p+2(p-1)}\eta_0\left(\frac{x-a_n}{r_n}\right)^{2\alpha p+2(p-1)} dx\\
& \sim r_n^{-2p} n^{2\alpha p+2(p-1)+\ell p}\int_{B_{r_n}(a_n)} \eta_0\left(\frac{x-a_n}{r_n}\right)^{2\alpha p+2(p-1)} dx\\
& = M_p r_n^{2-2p} n^{2\alpha p+2(p-1)+\ell p}\\
& = M_p e^{2(p-1)n^q}n^{2(\alpha + 1 -\ell)p - 2},
\end{align*}
with $M_p= \ds \int_{\R^2} \eta_0^p(z)dz >0$. We can verify easily \eqref{aK0} and $\alpha_1(K_0) = \frac{\ell-1}{2}$.

\medskip
The following result exhibits clearly the asymptotic behavior of Cheng-Lin's solutions to \eqref{eq 1.1} with $K= K_0$.
\begin{teo}\label{teo 4.1}
Assume that $\ell > q > 1$ and $K=K_0$ given by (\ref{k2}). By \cite{CL}, for any $0 < \alpha < \frac{\ell-1}{2}$, there is a unique solution of \eqref{eq 1.1} $u_\alpha$ satisfying \eqref{ChenL}. Denote $\alpha_* = \frac{\ell-q}{2}$, we have:
\begin{enumerate}
\item[(i)] For any $0 < \alpha < \alpha_*$, $u_\alpha$ satisfying \eqref{newb1}.
\item[(ii)] For $\alpha = \alpha_*$, the solution denoted by $u_*$ satisfies \eqref{b1}, but the remainder term $u_{\alpha_*}- \alpha_*\ln|x|$ does not converge to any constant as $|x|\to \infty$.
\item[(iii)] For any $\alpha_* < \alpha < \frac{\ell-1}{2}$,  the remainder term $u_{\alpha}  -\alpha\ln |x|$ is unbounded at infinity.
\end{enumerate}
\end{teo}

To our best knowledge, this is a first example of $K_0$ without any symmetry or growth assumption, where we show a complete picture of asymptotic behavior for solutions provided in \cite{CL}. This is also a first example of solution to \eqref{eq 1.1} with nonpositive Gaussian curvature, logarithmic growth but non uniform behavior at infinity.

\medskip
The rest of our paper is organized as follows.  In Section 2, we show the link between solutions satisfying \eqref{b1} and that provided by Cheng-Lin, we show also some properties of $\alpha_p(K)$ and the relationship between $({\mathbb H}_0)$ and $({\mathbb H}_1)$.  In Section 3, we prove Theorem \ref{teo 1}. Section 4 is addressed to the study of solutions with $K_0$. In the following, $C$ or $C_i$ denote always generic positive constants, their values could be changed from one line to another.

\setcounter{equation}{0}
\section{Preliminaries}
First, we indicate the relationship  between solutions satisfying \eqref{b1} and solutions given by \cite{CL}.
\begin{proposition}
\label{lnew1}
Let $u$ be a continuous solution to \eqref{eq 1.1} satisfying \eqref{b1} with $\alpha < \alpha_1(K)$. Then $u$ is the Cheng-Lin's solution, i.e. $u$ satisfies \eqref{ChenL}.
\end{proposition}

\noindent
{\bf Proof.} We involve function $w_0$ given by \eqref{defw}. So $\Delta w_0(x) =0$ for $|x|>1$ and
\begin{eqnarray}
\label{w0}
\int_{\R^2}\Delta  w_0 dx = 2\pi.
\end{eqnarray}
We claim that
 \begin{equation}\label{3.2.1}
\int_{\R^2} |\nabla(u - \alpha w_0)|^2 dx< \infty \aand  \int_{\R^2}  Ke^{2u}  dx = - 2\alpha \pi.
 \end{equation}

Denote $v = u - \alpha w_0$, we prove first $v \in L^2(\R^2)$. To this end, set $\eta_\rho(t)=\eta_0(t/\rho)$ for $\rho>0$,
where $\eta_0$ is a smooth function satisfying \eqref{eta 0}. Observe that
\begin{equation}\label{eq 3.1}
-\Delta v = K e^{2\alpha w_0} e^{2v} + \alpha \Delta  w_0 \quad{\rm in}\;\; \mathbb{R}^2.
\end{equation}
As $v$ is bounded in $\R^2$ and $\alpha < \alpha_1(K)$, there holds
\begin{align*}
2\int_{\R^2} \eta_\rho^2|\nabla v|^2 dx & = \int_{\R^2} \eta_\rho^2\left(\Delta v^2 -2v\Delta v \right) dx\\
& = \int_{\R^2} v^2\Delta \eta_\rho^2dx- 2\int_{\R^2} \eta_\rho^2v\Delta v dx\\
& \leq \|v\|_\infty^2 \int_{\R^2} |\Delta \eta_\rho^2| dx + 2 \|v\|_\infty\int_{\R^2} |\Delta v|dx\\
 & \leq C_1 + 2 \|v\|_\infty \int_{\R^2} \left(|K|e^{2\alpha w_0+2\|v\|_\infty} + |\alpha\Delta  w_0|\right) dx \leq C_2.
\end{align*}
Here $C_i$ are positive constants independent on $\rho$. Passing $\rho\to \infty$, we see that $\nabla v \in L^2(\R^2)$.

\medskip
Furthermore,
\begin{align*}
\left|\int_{\R^2} \eta_\rho \Delta v dx\right| & = \left|\int_{\R^2} \nabla \eta_\rho\cdot\nabla v dx\right|\\
 &\le  \left(\int_{B_{2\rho}\setminus B_{\rho}}  |\nabla v|^2 dx\right)^{\frac12}  \left(\int_{B_{2\rho}\setminus B_{\rho}}  |\nabla \eta_\rho |^2 dx\right)^{\frac12}
 \\& =  \|\nabla\eta_0\|_{L^2(\R^2)}\left(\int_{B_{2\rho}\setminus B_{\rho}}  |\nabla v|^2 dx\right)^{\frac12}.
\end{align*}
Remark that $\Delta v \in L^1(\R^2)$ and $\nabla v \in L^2(\R^2)$, we obtain, taking $\rho\to \infty$,
$$\int_{\R^2} \Delta v dx =0.$$
Combining \eqref{eq 3.1} and \eqref{w0}, we obtain the second part of (\ref{3.2.1}). \hfill$\Box$

\medskip
We give some elementary consequences of $(\mathbb{H}_1)$ and $(\mathbb{H}_0)$.
\begin{lemma}\label{pr 2.1}
Assume that $(\mathbb{H}_1)$ holds true with $p_0 > 1$. Then for any $p\in[1,p_0)$, $\alpha_p(K)\ge \alpha_{p_0}(K)$ and $\alpha_p(K)\to \alpha_1(K)$ as $p\to 1^+$.
\end{lemma}

\noindent
{\bf Proof.} For simplicity, we denote $\alpha_p(K)$ by $\alpha_p$. Take any $\beta < \alpha<\alpha_{p_0}$ and $p\in [1,\,p_0)$. Denote $\tau_\alpha=-2+\frac{2p_0p}{p_0-p}(\beta-\alpha)<-2$, then
$$\tau_\alpha\left(1-\frac{p}{p_0}\right) + \Big[2\alpha p_0+2(p_0-1)\Big]\frac{p}{p_0}=2\beta p+2(p-1).$$
By H\"older's inequality, as $\alpha < \alpha_{p_0}$,
\begin{align*}
&\int_{\R^2}  |K(x)|^p(1+|x|)^{2\beta p+2(p-1)} dx\\
\le& \; \left(\int_{\R^2}   |K(x)|^{p_0}(1+|x|)^{2\alpha p_0+2(p_0-1)}  dx\right)^{\frac{p}{p_0}}\left(\int_{\R^2}(1+|x|)^{\tau_\alpha} dx\right)^{1-\frac{p}{p_0}} < \infty.
\end{align*}
It means that $\alpha_p\ge \alpha_{p_0}$ for any $p\in [1,\,p_0)$. We see also that the function $p \mapsto \alpha_p$ is non-increasing in $[1, p_0]$.

\medskip
If $\alpha_{p_0}=\alpha_1$, then $\alpha_p \equiv \alpha_{p_0}$ for $p\in[1,p_0]$ and we are done. Now assume that $\alpha_{p_0}<\alpha_1$. Let $\e \in (0, \alpha_1 - \alpha_{p_0})$, take $$\l = \alpha_1 - \e, \quad \beta_2<\alpha_{p_0} <\l <\beta_1<\alpha_1,$$
and consider
\begin{equation*}
h(p) :=  \beta_1(p_0-p) + \beta_2 p_0(p-1) - \l p(p_0-1)\quad \mbox{in } \;[1, p_0].
\end{equation*}
Readily, $h(1) =(\beta_1-\l)(p_0-1) > 0$ and $h(p_0) = (\beta_2-\l)(p_0-1)p_0<0.$ Hence, there exists $p_* \in (1, p_0)$ such that $h(p_*) = 0$, or equivalently
\begin{align*}
 2\beta_1\theta + \Big[2\beta_2 p_0+2(p_0-1)\Big](1-\theta)=2\l p_*+2(p_*-1) \quad \mbox{where }\; \theta :=\frac{p_0-p_*}{p_0-1} \in (0, 1).
\end{align*}
As $\theta+p_0(1-\theta)=p_*$, there holds
\begin{align*}
 &\int_{\R^2}  |K(x)|^{p_*}(1+|x|)^{2\l p_*+2(p_*-1)}  dx\\
\le & \; \left(\int_{\R^2}  |K(x)|(1+|x|)^{2\beta_1 } dx\right)^{\theta} \left(\int_{\R^2} |K(x)|^{p_0}(1+|x|)^{2\beta_2 p_0+2(p_0-1)}dx\right)^{1-\theta} < \infty.
\end{align*}
Thus we get $\alpha_{p_*} \geq \l$, so $\alpha_1 - \e \leq \alpha_p \leq \alpha_1$ for $p\in [1, p_*]$, which means $\lim_{p\to 1^+} \alpha_p = \alpha_1$. \hfill$\Box$

\begin{remark}
By Lemmas \ref{pr 2.1}, we see that $({\mathbb H}_1)$ is equivalent to $({\mathbb H}_1')$: There exists $p > 1$ such that $\alpha_p(K) > 0$.
\end{remark}

\begin{lemma}\label{pr 2.3}
If $K$ satisfies $({\mathbb H}_0)$, then there exists $p>1$ such that $\alpha_p(K) > 0$.
\end{lemma}

\noindent
{\bf Proof.} Fix some $\alpha\in(0,\alpha_1)$ and choose $p>1$ but nearby such that
$$\beta:=\alpha+ \left(\frac m2+\alpha+1\right)(p-1)<\alpha_1.$$
Then
\begin{align*}
|K(x)|^{p-1}(1+|x|)^{2(\alpha+1)(p-1) } \le C(1+|x|)^{m(p-1)+2(\alpha+1)(p-1)} = C(1+|x|)^{2\beta-2\alpha},
\end{align*}
so that
\begin{align*}
\int_{\R^2}  |K(x)|^p(1+|x|)^{2\alpha p+2(p-1)} dx & = \int_{\R^2}  |K(x)| (1+|x|)^{2\alpha} \Big[|K(x)|^{p-1}(1+|x|)^{2(\alpha+1)(p-1)}\Big] dx   \\
 & \le C\int_{\R^2}  |K(x)| (1+|x|)^{2\beta} dx < \infty.
\end{align*}
The proof is completed.\hfill$\Box$

\medskip
We end this section by a simple fact of uniqueness, which can be derived from Theorem 4.1 in \cite{CN1}. For the convenience of readers, we give a short proof here.
\begin{lemma}
Let $K \leq 0$, then for any $\alpha \in \R$, there exists at most one continuous solution to \eqref{eq 1.1} satisfying \eqref{b1}.
\end{lemma}
Let $\widetilde u_\alpha$ and $u_\alpha$ be two solutions of \eqref{eq 1.1} satisfying \eqref{b1} with $K \leq 0$.
Denote $v=\widetilde u_\alpha-u_\alpha$, then
\begin{align*}
 -\Delta v^2 = -2v\Delta v - 2|\nabla v|^2\le -2v\Delta v   = 2Ke^{2u_\alpha}v \left(e^{2v}-1\right)
 \le 0,
\end{align*}
that is, $v^2$ is bounded and subharmonic in $\R^2$, hence $v^2$ is a constant. So $0 = -\Delta v = Ke^{2u_\alpha}(e^{2v}-1)$ in $\R^2$, which implies that $v \equiv 0$ and the uniqueness follows.

\section{Proof of Theorem \ref{teo 1} }
\setcounter{equation}{0}
We begin with a result which plays an important role to obtaining the asymptotic behavior (\ref{1.3}).

 \begin{proposition}\label{lm 3.1}
Let
\begin{align*}
f \in \dot{L}^1(\R^2) := \left\{ f \in L^1(\R^2), \; \int_{\R^2} f(x) dx = 0 \right\}.
\end{align*}
Assume that $\alpha_{p_0}(f) > 0$ with $p_0 > 1$, then there exists a unique ground state solution to $-\Delta w = f$ in $\R^2$. Moreover, $w \in C(\R^2)$ and
\begin{equation}\label{4.2}
\lim_{|x|\to \infty} w(x)|x|^{\frac{2\beta}{1+2\beta}} = 0 \quad {\rm for\ any}\; \beta\in(0,\alpha_1(f)).
\end{equation}
\end{proposition}
\noindent {\bf Proof.} The uniqueness of the ground state solution is a trivial consequence of Liouville's Theorem.
So we prove only the existence of a solution satisfying \eqref{4.2}. By $\alpha_{p_0}(f)>0$ and Lemma \ref{pr 2.1},
we get $\alpha_1(f)>0$ and $f(z)\ln(2+|z|) \in L^1(\R^2)$. As $f \in \dot{L}^1\cap L^{p_0}_{loc}$, consider
\begin{align*}
w(x) = -\frac{1}{2\pi}\int_{\R^2}\ln|x - y|f(y) dy \quad  \mbox{in } \R^2.
\end{align*}
Obviously, $-\Delta w = f$ in $\R^2$, $w \in C(\R^2)$ and for any $x \ne 0$,
\begin{align*}
2\pi w(x) = \frac{1}{2\pi}\int_{\R^2}\ln\frac{|x|}{|x - y|}f(y) dy = \sum_{i=1}^3 \int_{\Omega_i} \ln\frac{|x|}{|x - y|}f(y) dy  =:  \sum_{i=1}^3 J_i(x),
\end{align*}
where $\Omega_1 = B_R$, $R \ge 2$, $\Omega_2 = B_{\frac{|x|}{2}}(x)$ and $\Omega_3 = \R^2\backslash(\Omega_1\cup\Omega_2)$. For $|x| \ge 4R > 0,$
\begin{align*}
|J_1(x)| &\le C\int_{B_R} \frac{|y|}{|x|} |f(y)|dy \le \frac{CR}{|x|} \norm{f}_{L^1(\R^2)}.
\end{align*}
Let $e_x = \frac{x}{|x|}$, $p_0'=\frac{p_0}{p_0-1}$ and fix $0 < \beta < \alpha_{p_0}(f)$, we get
\begin{align*}
|J_2(x)| &=  \left|\int_{\Omega_2}  \ln\frac{|x|}{|x - y|}f(y) dy\right|
\\&\le
\left(\int_{\Omega_2} |f(y)|^{p_0} (1+|y|)^{2\beta p_0+2(p_0-1)}  dy\right)^{\frac1{p_0}} \left(\int_{\Omega_2} \left|\ln\frac{|x|}{|x - y|}\right|^{p_0'} (1+|y|)^{-2- 2\beta p_0'} dy\right)^{\frac1{p_0'}}
\\&\le C\left(\int_{B_{\frac12}(e_x)}\big|\ln|e_x - z|\big|^{p_0'} |z|^{-2-2\beta p_0'} |x|^{-2\beta p_0'} dz\right)^{\frac1{p_0'}}
\\& = C|x|^{-2\beta}.\end{align*}
Let $|x| \ge 4R \geq 8$. In $\Omega_3$, as $\min(|x - y|, |x|, |y|) \ge R \ge 2$, we can check that
$$\forall\; y \in \Omega_3, \qquad \left|\ln\frac{|x|}{|x - y|}\right| \le \ln|y|,$$
since $\frac{1}{|y|} \leq \frac{|x|}{|x - y|} \leq |y|$.
Then there holds
\begin{align*}
|J_3(x)| & \le \int_{\Omega_3} |f(y)|\ln |y| dy\\
& \le  \left(\int_{\Omega_3} |f(y)|^{p_0}(1+|y|)^{2\beta p_0+2(p_0-1)}   dy\right)^{\frac1{p_0}}\left(\int_{\Omega_3} (\ln|y|)^{p_0'}  |y|^{-2-2\beta p_0'} dy\right)^{\frac1{p_0'}}\\
& \le C\left(\int_{\R^2\backslash B_R} (\ln|y|)^{p_0'}  |y|^{-2-2\beta p_0'} dy\right)^{\frac1{p_0'}}\\
& \le C R^{-2\beta}\ln R.
\end{align*}
Notice that all the constants $C$ involved are independent of $|x| \geq 4R$ and $R\geq 2$. Taking $R= |x|^{\frac1{1+2\beta }}$, there holds
 \begin{equation}\label{4.4-0}
  |w(x)|\le C|x|^{-\frac{2\beta }{1+2\beta }}\ln|x|\quad {\rm for } \ |x| \; {\rm large}.
 \end{equation}
This ends the proof of \eqref{4.2}, because we can get \eqref{4.4-0} with any $p \in (1, p_0)$ and any $\beta < \alpha_p(f)$. \hfill$\Box$

\bigskip
We are now in position to prove Theorem \ref{teo 1}.

\medskip
\noindent {\bf Proof of Theorem \ref{teo 1}.} Let $\alpha\in(0,\alpha_1(K))$ and $w_0$ be given by \eqref{defw},
 there exists then a unique $t_\alpha \in \R$ such that
$$g_\alpha := K e^{2\alpha w_0} e^{2t_\alpha} + \alpha \Delta  w_0 \in \dot{L}^1(\R^2).
$$
Indeed, as $K \leq 0$ in $\R^2$ and $\alpha > 0$, the unique choice is given by
$$t_\alpha = \frac{1}{2}\ln\left(-\frac{2\alpha\pi}{\ds \int_{\R^2} K e^{2\alpha w_0} dx}\right).
$$
On the other hand, using $(\mathbb{H}_1)$ and Lemma \ref{pr 2.1}, there exists $p > 1$ such that $\alpha < \alpha_p(K)$. As $\alpha\Delta w_0$ is compactly supported, it's not difficult to see that
$$\alpha_p(g_\alpha) = \alpha_p(K) - \alpha > 0.$$
Using Proposition \ref{lm 3.1}, we get a continuous ground state solution to $-\Delta \widetilde w = g_\alpha$ in $\R^2$.

\medskip
We conclude by super-sub solution method. Let $\overline v = \alpha w_0 + t_\alpha + \widetilde w + \|\widetilde w\|_\infty$. As $K \leq 0$, it's easy to see that
$$-\Delta \overline v = -\alpha\Delta w_0 - \Delta \widetilde w= K e^{2(\alpha w_0+t_\alpha)} \geq K e^{2\overline v}.$$
Hence $\overline v$ is a supersolution to equation \eqref{eq 1.1}. Similarly,
$\underline v = \alpha w_0 + t_\alpha + \widetilde w - \|\widetilde w\|_\infty$
is a subsolution to \eqref{eq 1.1}. Obviously $\overline v \geq \underline v$ in $\R^2$.
Using the standard Perron's method, there exists a solution $u$ to \eqref{eq 1.1}
such that $\underline v \le u \le \overline v$ in $\R^2$.  Clearly, $u$ satisfies \eqref{b1}.

\medskip
By Proposition \ref{lnew1}, there holds $Ke^{2u} + \alpha\Delta w_0  \in \dot{L}^1(\R^2)$. Applying Lemma \ref{pr 2.1}, there exists $p > 1$ such that
$$\alpha_p\left(Ke^{2u} + \alpha\Delta w_0\right) = \alpha_p(K) - \alpha > 0.$$
By Proposition \ref{lm 3.1}, we get a continuous solution of $-\Delta w = Ke^{2u} + \alpha\Delta w_0$ in $\R^2$,
which satisfies \eqref{4.2} for any $0 < \beta < \alpha_1(K) - \alpha$.
Therefore, $v = u - \alpha w_0 - w$ is a bounded harmonic function over $\R^2$, hence $v$ is a constant $c_\alpha$.
Finally $u = \alpha w_0 + c_\alpha + w$, so \eqref{1.3} holds trues for any $0 < \beta < \alpha_1(K) - \alpha$. \hfill$\Box$

\medskip
Similarly, using the super-sub solution method and Proposition \ref{lnew1}, we can claim
\begin{proposition}
\label{teo 2}
Assume that for $K\le 0$, $\alpha > 0$ and $t \in \R$, there exists a bounded solution to $-\Delta v = Ke^{2\alpha w_0 + 2t} + \alpha\Delta w_0$ in $\R^2$, then the equation (\ref{eq 1.1}) admits a unique solution $u_\alpha$ satisfying \eqref{b1}. Moreover, the total curvature of $u_\alpha$ is $-2\alpha\pi$.
\end{proposition}

\setcounter{equation}{0}
\section{Exotic curvature $K_0$ and anisotropic behavior solutions}
Here we prove Theorem \ref{teo 4.1}. We divide our study by several steps.

\subsection{Solution satisfying \eqref{b1} for $\alpha \leq \alpha_*$.}
\label{c1}
Remark that we cannot apply Theorem 1.1 in \cite{CL} or its proof, since they don't provide the asymptotic behavior. Indeed, we will construct directly solutions satisfying \eqref{b1}, and Proposition \ref{lnew1} ensures then they are solutions given by \cite{CL}. Remark also that we cannot apply the results of Section 3, since $\alpha_p(K_0) = -\infty$ for any $p > 1$.

\medskip
Recall that $\alpha_1(K_0) = \frac{\ell - 1}2 > 0$ and $w_0$ is given by \eqref{defw}. Let $0 < \alpha \le \frac{\ell-q}{2}$ and
\begin{align*}
v_\alpha(x) = -\frac{1}{2\pi}\int_{\R^2} K_0(y) e^{2\alpha w_0(y)}\ln|x-y|dy.
\end{align*}
It is easy to see that $v_\alpha$ is well defined, continuous and locally bounded in $\R^2$. Let
\begin{align}
\label{new4}
n \geq 2 \aand x\in B_{n+\frac12}\setminus \overline B_{n-\frac12}.
\end{align}
We decompose $v_\alpha(x)$ as follows,
\begin{align*}
2\pi v_\alpha(x)  & = \sum^{\infty}_{k=2,\,k\not=n}r_k^{-2}k^{-\ell} \int_{B_{r_k}(a_k)}  |y|^{2\alpha} \eta_0\left(\frac{y - a_k}{r_k}\right)\ln|x-y| dy
\\& \hspace*{0.5cm} + r_n^{-2}n^{-\ell} \int_{B_{r_n}(a_n) } |y|^{2\alpha} \eta_0\left(\frac{y - a_n}{r_n}\right)\ln|x-y|\,dy\\
& =: E_1(x)+E_2(x).
\end{align*}

Denote
\begin{align*}
s_k :=\int_{B_{r_k}(a_k)} |y| ^{2\alpha} \eta_0\left(\frac{y - a_k}{r_k}\right) dy \le C(k+1)^{2\alpha} r_k^2, \quad \forall\; k \geq 2.
\end{align*}
Consider first $k\le\frac n2$ and $k \ge 2$. Notice that there exists $C > 0$ such that
$$\Big|\ln|x-y|-\ln|x|\Big|\le C\frac{|y|}{|x|}, \qquad \forall\; |x|\ge 2|y| > 0.$$
For $x$ satisfying \eqref{new4}, there holds then
\begin{align*}
\left|\int_{B_{r_k}(a_k) } |y| ^{2\alpha} \eta_0\left(\frac{y - a_k}{r_k}\right) \ln|x-y|\,dy- s_k \ln|x| \right|
& = \left|\int_{B_{r_k}(a_k) }  |y| ^{2\alpha} \eta_0\left(\frac{y - a_k}{r_k}\right) \ln\frac{|x-y|}{|x|} dy\right|\\
& \le \frac{C}{n} (k+1)^{2\alpha+1} \int_{B_{r_k}(a_k) } \eta_0\left(\frac{y - a_k}{r_k}\right)  dy
\\  & = \frac{C}{n}  k^{2\alpha+1} r_k^2.
\end{align*}

Let now $k\ne n$, $k > \frac{n}{2}$. As $|y|\geq \frac{7}{4}$ and $|x - y| \ge \frac{1}{4}$ in $B_{r_k}(a_k)$, there exists $C > 1$ (independent of $k, n\ge 2$, $k \ne n$) such that
$$\forall\; y \in B_{r_k}(a_k), \quad |y|^{-C}\leq \frac{\frac{1}{4}}{|y| + \frac{1}{4}}\le \frac{|x-y|}{|y| + |x -y|} \le \frac{|x-y|}{|x|} \le 1+ |y| \le |y|^C.$$
We obtain
\begin{align*}
\left|\int_{B_{r_k}(a_k) } |y| ^{2\alpha} \eta_0\left(\frac{y - a_k}{r_k}\right) \ln|x-y|\,dy- s_k \ln|x| \right| & = \left|\int_{B_{r_k}(a_k) }  |y| ^{2\alpha} \eta_0\left(\frac{y - a_k}{r_k}\right) \ln\frac{|x-y|}{|x|} dy\right|
\\ &\le  C\left|\int_{B_{r_k}(a_k) }  |y|^{2\alpha}  \eta_0\left(\frac{y - a_k}{r_k}\right)\ln |y| dy\right|  \\
   &\le   C (k+1)^{2\alpha} r_k^2 \ln (k+1).
\end{align*}
Therefore, for $x$ in \eqref{new4},
\begin{align*}
& \left|E_1(x) - \left(\sum_{k=2}^{\infty} r_k^{-2}k^{-\ell} s_k\right)  \ln|x| \right| \\
\le & \; \frac{C}{n} \sum^{[\frac n2]}_{k=2}k^{2\alpha-\ell +1}
+ C \sum_{k > [\frac n2], k \ne n}k^{2\alpha-\ell} \ln k + r_n^{-2}n^{-\ell} s_n\ln|x|\\
\le & \;  \frac{C}{n} \sum^{[\frac n2]}_{k=2}k^{2\alpha-\ell +1}
+ C \sum_{k = [\frac n2]+1}^\infty k^{2\alpha-\ell} \ln k \quad \rightarrow 0, \;\;\mbox{as } n\to \infty.
 \end{align*}
Here $[\frac n2]$ is the integer part of $\frac n2$. Finally, as $2\alpha  - \ell + 1 \leq 1 - q < 0$, there holds
\begin{align}
\label{E1}
\lim_{|x| \to \infty} \Big|E_1(x) - A_\alpha \ln|x| \Big| = 0, \quad \mbox{where } \; A_\alpha = \sum_{k=2}^{\infty} r_k^{-2}k^{-\ell} s_k \in (0, \infty).
\end{align}

On the other hand, for $x$ satisfying \eqref{new4},
\begin{align}
\label{E2}
\begin{split}
|E_2(x)|&\le  Cr_n^{-2}n^{2\alpha -\ell}  \int_{B_{r_n}(a_n) }  \eta_0\left(\frac{y - a_n}{r_n}\right) \times (-\ln|x-y|) dy\\
& = -Cn^{2\alpha -\ell}\int_{B_1} \eta_0(z)\ln|x-a_n-r_nz| dz\\
&\le -Cn^{2\alpha -\ell}\int_{B_1} \eta_0(z)\ln|r_n z| dz\\
&=  Cn^{2\alpha -\ell}(|\ln r_n| + 1)\\
& \le Cn^{2\alpha -\ell+q}.
\end{split}
\end{align}
For the third line, we used the symmetric decreasing rearrangement argument, as $\eta_0$ and $-\ln|r_n z|$ are radial non-increasing positive functions in $B_1$.

\medskip
Combining \eqref{E1}--\eqref{E2}, $v_\alpha$ is a solution to $-\Delta v_\alpha = K_0e^{2\alpha w_0}$ in $\R^2$ such that
\begin{align*}
v_\alpha(x) = \frac{A_\alpha}{2\pi} \ln|x| + O(1), \; \; \mbox{as } |x| \to \infty.
\end{align*}
Taking
$$t_\alpha =\frac12 \ln\frac{2\alpha\pi}{A_\alpha},$$
$h_\alpha = e^{2t_\alpha}v_\alpha - \alpha w_0$ is clearly a bounded solution to
\begin{align*}
-\Delta h = K_0 e^{2\alpha w_0 + 2t_\alpha} + \alpha\Delta w_0 \quad \mbox{in } \; \R^2.
\end{align*}
By Propositions \ref{teo 2}, for any $\alpha\in(0,\, \frac{\ell-q}{2}]$, equation (\ref{eq 1.1}) with $K = K_0$ admits a unique solution satisfying \eqref{b1}, and whose total curvature is equal to $-2\alpha\pi$.

\subsection{Uniform behavior for $0 < \alpha < \alpha_*$.}
For $\alpha \in (0, \frac{\ell - q}{2})$, we have a continuous solution $u_\alpha$ to \eqref{eq 1.1} satisfying \eqref{b1}. Rewrite $u_\alpha = \alpha w_0 + \xi_\alpha$ with $\xi_\alpha \in L^\infty(\R^2)$. Let
\begin{align*}
\widetilde v_\alpha(x) = -\frac{1}{2\pi}\int_{\R^2} K_0(y) e^{2\alpha w_0(y)+ 2\xi_\alpha(y)}\ln|x-y|dy.
\end{align*}
By decomposing $\widetilde v_\alpha$ as above and using
\begin{align*}
\widetilde s_k :=\int_{B_{r_k}(a_k)} e^{2\xi_\alpha(y)}|y| ^{2\alpha} \eta_0\left(\frac{y - a_k}{r_k}\right) dy \quad \mbox{for }\; k \geq 2,
\end{align*}
the similar estimates to \eqref{E1}--\eqref{E2} hold true for $\widetilde v_\alpha$, which yield that (as now $2\alpha -\ell+q < 0$)
\begin{align}
\label{wv}
\lim_{|x| \to \infty}\left[\widetilde v_\alpha(x) - \frac{\widetilde A_\alpha}{2\pi} \ln|x| \right] = 0, \quad \mbox{where }\; \widetilde A_\alpha =  \sum_{k=2}^{\infty} r_k^{-2}k^{-\ell} \widetilde s_k \in (0, +\infty).
\end{align}
So $u_\alpha - \widetilde v_\alpha$ is a harmonic function in $\R^2$ and $(u_\alpha- \widetilde v_\alpha)(x) = O(\ln|x|)$ at infinity, hence $u_\alpha - \widetilde v_\alpha$ is a constant in $\R^2$. The behavior \eqref{wv} means that $u_\alpha$ satisfies \eqref{newb1}.

\subsection{Anisotropic behavior for $\alpha = \alpha_*$.}
Let $\alpha = \alpha_* = \frac{\ell - q}{2}$, and $u_*$ be the unique solution to \eqref{eq 1.1} satisfying \eqref{b1}, hence $\xi_* := u_* - \alpha_*w_0$ is uniformly bounded in $\R^2$. However, we will show that
$\xi_*$  does not converge to a constant at infinity. Define
\begin{align*}
2\pi v_*(x) &= -\int_{\R^2} K_0(y) e^{2\alpha_* w_0(y) + 2\xi_*(y)}\ln|x-y|dy \\
& = \sum^{\infty}_{k=2,\,k\not=n}r_k^{-2}k^{-\ell} \int_{B_{r_k}(a_k)} e^{2\xi_*(y)}  |y|^{2\alpha} \eta_0\left(\frac{y - a_k}{r_k}\right)\ln|x-y|\,dy
\\&\hspace*{0.5cm} + r_n^{-2}n^{-\ell} \int_{B_{r_n}(a_n) } e^{2\xi_*(y)} |y|^{2\alpha} \eta_0\left(\frac{y - a_n}{r_n}\right)\ln|x-y| dy
\\&=: F_1(x)+F_2(x).
\end{align*}
We see that
\begin{align*}
-\Delta v_* = K_0e^{2\alpha_* w_0+ 2\xi_*} = K_0 e^{2u_*}\;\; \mbox{in } \; \R^2 \aand v_*(x) = O(\ln|x|) \; \; \mbox{as }\; |x|\to\infty.
\end{align*}
The logarithmic control of $v_*$ at infinity can be obtained as above, using $\xi_* \in L^\infty(\R^2)$.

\medskip
We shall estimate the difference between $v_*(a_n)$ and $v_*(-a_n)$. The estimation of $F_1$ is very similar to $E_1$ in the above. Let $n \geq 2$. If $k\le\frac n2$ and $k \ge 2$, there holds
\begin{align*}
 &\left|\int_{B_{r_k}(a_k) } e^{2\xi_*(y)} |y|^{2\alpha_*} \eta_0\left(\frac{y - a_k}{r_k}\right) \Big[\ln|a_n - y|-\ln|a_n+ y| \Big] dy \right| \\
  \le &\;\frac{C}{n}  e^{2\norm{\xi_*}_\infty}(k+1)^{2\alpha_*} \int_{B_{r_k}(a_k) } \eta_0\left(\frac{y - a_k}{r_k}\right) |y| \,dy
\\ \le& \; \frac{C}{n}  k^{2\alpha_*+1} r_k^2.
\end{align*}
Let $k>\frac n2$ and $k\not=n$, we have
\begin{align*}
 & \left|\int_{B_{r_k}(a_k) } e^{2\xi_*(y)}  |y|^{2\alpha_*} \eta_0\left(\frac{y - a_k}{r_k}\right)\Big[\ln|a_n - y|-\ln|a_n+ y| \Big] dy \right| \\
 \le  & \; C(k+1)^{2\alpha_*} (\ln n+\ln k )\int_{B_{r_k}(a_k) } \eta_0\left(\frac{y - a_k}{r_k}\right)  \,dy\\
\le & \; Ck^{2\alpha_*} r_k^2 \ln k.
\end{align*}
As $2\alpha_* - \ell = -q < -1$, we get then
\begin{align}
\label{new2}
|F_1(a_n)-F_1(-a_n)|&\le \frac{C}{n}  \sum^{[\frac n2]}_{k=2}k^{2\alpha_*-\ell +1}
+ C\sum^{\infty }_{k=[\frac n2]+1}k^{2\alpha_*-\ell} \ln k \;\; \rightarrow 0,\quad {\rm as}\;\; n\to\infty.
\end{align}

\medskip
Now we estimate $F_2(\pm a_n)$:
\begin{align*}
 0\le F_2(-a_n)  &\le  Cr_n^{-2} n^{-\ell+2\alpha_* } \int_{B_{r_n}(a_n) }   \eta_0\left(\frac{y-a_n}{r_n}\right)  \ln|a_n+ y| dy \le  Cn^{-\ell+2\alpha} \ln (2n+1)
\end{align*}
and
\begin{align*}
 -F_2(a_n)  & = -r_n^{-2}n^{-\ell} \int_{B_{r_n}(a_n) } e^{2\xi_*(y)} |y|^{2\alpha_*} \eta_0\left(\frac{y - a_n}{r_n}\right)\ln|a_n - y| dy\\
& \ge  -e^{-2\norm{\xi_*}_\infty} r_n^{-2}n^{-\ell}(n - 1)^{2\alpha_*}  \int_{B_{r_n}(a_n)}   \eta_0\left(\frac{y - a_n}{r_n}\right)  \ln|a_n - y| dy
\\ &\ge  Cn^{-\ell+2\alpha_*} \int_{B_1}   \eta_0(z)\Big[-\ln r_n - \ln|z|\Big]dz
\\& \ge -Cn^{-\ell+2\alpha_*}\ln r_n\\
& = Cn^{-\ell+2\alpha_*+ q} = C > 0.
\end{align*}
Therefore $\liminf_{n\to \infty}|F_2(a_n) - F_2(-a_n)| \geq C > 0$. Combining with \eqref{new2}, we get
\begin{align*}
\liminf_{n\to \infty}|v_*(a_n) -v_*(-a_n)| \geq \frac{C}{2\pi} > 0.
\end{align*}
Furthermore, as $-\Delta (v_* - u_*)=0$ in $\R^2$ and $v_* - u_*$ is $O(\ln |x|)$ at infinity, the standard Liouville's theorem implies that $v_* - u_*$ is constant in $\R^2$. Consequently,
\begin{align*}
\liminf_{n\to \infty}|\xi_*(a_n) - \xi_*(-a_n)|  = \liminf_{n\to \infty}|v_*(a_n) -v_*(-a_n)| \geq \frac{C}{2\pi} > 0.
\end{align*}
The above estimate means that $\xi_*$ does not converge to a constant at infinity.

\subsection{Unbounded remainder term for $\alpha_* < \alpha < \alpha_1(K_0)$.}  Recall that by Theorem 1.1 in \cite{CL}, for any $\alpha\in \left(\frac{\ell-q}{2},\, \frac{\ell-1}{2}\right)$,  problem (\ref{eq 1.1}) possesses a unique solution $u_\alpha$ satisfying \eqref{ChenL}. We shall prove that $\xi_\alpha := u_{\alpha} -\alpha w_0$ is unbounded.

\medskip
Suppose the contrary, assume that $\xi_\alpha$ is bounded in $\R^2$. Let
$$\eta_\alpha(x) = -\frac{1}{2\pi}\int_{\R^2} K_0(y) e^{2u_\alpha(y)} \ln|x-y|dy \quad \mbox{in } \; \R^2,$$
then we have
\begin{align*}
2\pi \eta_\alpha (x) &= \sum^{\infty}_{k=2,\,k\not=n}r_k^{-2}k^{-\ell} \int_{B_{r_k}(a_k) }  |y|^{2\alpha} e^{2\xi_\alpha} \eta_0\left(\frac{y - a_k}{r_k}\right)\ln|x-y| dy
\\& \hspace*{0.5cm} + r_n^2n^{-\ell} \int_{B_{r_n}(a_n)}  |y|^{2\alpha} e^{2\xi_\alpha} \eta_0\left(\frac{y - a_n}{r_n}\right)\ln|x-y| dy
\\& =: G_1(x) +G_2(x).
\end{align*}
We will estimate $\eta_\alpha(a_n)$. For $k\ne n$, $k \ge 2$,
 \begin{align*}
&\left|\int_{B_{r_k}(a_k) } |y| ^{2\alpha}e^{2\xi_\alpha} \eta_0\left(\frac{y - a_k}{r_k}\right) \ln|a_n - y|\,dy\right|\\
\le & \;  Ce^{2\|\xi_\alpha\|_\infty} (k+1)^{2\alpha} \int_{B_{r_k}(a_k) } \eta_0\left(\frac{y - a_k}{r_k}\right) \big(\ln n+\ln|y|\big) dy\\
\le & \;  C k^{2\alpha} r_k^2 \big(\ln n +\ln k\big).
\end{align*}
Since $\alpha < \frac{\ell-1}{2}$, there exists $C > 0$ such that $|G_1(a_n)| \le C\ln n$ for $n$ large. On the other hand,
\begin{align*}
-G_2(a_n)  &\ge -e^{-2\|\xi_\alpha\|_\infty} r_n^{-2}n^{-\ell}(n-1)^{2\alpha } \int_{B_{r_n}(a_n)}   \eta_0\left(\frac{y - a_n}{r_n}\right) \ln|a_n - y| dy
\\ &\ge  -Cn^{-\ell+2\alpha} \ln r_n\\
& = C_3 n^{-\ell+2\alpha+q}.
\end{align*}
Combining the estimations of $G_1(a_n)$ and $G_2(a_n)$, we get, as $\alpha > \frac{\ell - q}{2}$,
\begin{align}
\label{new3}
\lim_{n\to \infty}\frac{\eta_\alpha(a_n)}{\ln n} = -\infty.
\end{align}

Moreover, remark that for $x$ satisfying \eqref{new4} and $y \in B_{r_n}(a_n)$, there holds $|x - y| \leq 3|x|$. Similar estimations as above yield that there exists $C>0$ such that
$$\eta_\alpha(x) \le C\ln(2+ |x|) \;\; \mbox{in }\; \R^2.$$
As $\widetilde v := \eta_\alpha - u_\alpha = \eta_\alpha - \alpha w_0 - \xi_\alpha$ is harmonic over $\R^2$, and $\widetilde v(x) \le C\ln(2+|x|)$ in $\R^2$, $\widetilde v$ is constant by Liouville's theorem. However, $\eta_\alpha = \alpha w_0 + \xi_\alpha + C$ with $\xi_\alpha \in L^\infty$ contradicts the estimate \eqref{new3}. This means that $\xi_\alpha$ cannot be bounded in $\R^2$. \hfill$\Box$\medskip

\bigskip
\noindent{\small {\bf Acknowledgements.} H.C. is supported by NSFC (No.~11726614 and 11661045),
and Jiangxi Provincial Natural Science Foundation (No. 20161ACB20007).
F.Z. and D.Y. are supported by Science and Technology Commission of Shanghai Municipality (STCSM), grant No. 18dz2271000.
F.Z. is also supported by NSFC (No.~11726613 and 11431005).


\begin{thebibliography}{99}

\bibitem{Ah} L.V. Ahlifors,  An extension of Schwartz's lemma, {\it Trans. Amer. Math. Soc. 43}, 359-364 (1938).

\bibitem {CL} K.-S. Cheng and C.-S. Lin, Conformal metrics   with prescribed nonpositive Gaussian on $\R^2$,
 {\it  Calc. Var. PDE   11},  203-231 (2000).

%\bibitem{CLin}  K.-S. Cheng  and J.-T. Lin, On the elliptic equations $\Delta u = K(x)u^¦Ò$ and $\Delta u = K(x)e^{2u}$, {\it Trans. Amer. Math. Soc. 304}, 639-667 (1987).

\bibitem {CN1} K.-S. Cheng and W.-M. Ni, On the structure of the conformal Gaussian curvature equation on $\R^2$, {\it Duke Math. J. 62}, 721-737 (1991).

\bibitem {CN2} K.-S. Cheng and W.-M. Ni,  On the structure of the conformal Gaussian curvature equation on $\R^2$ II, {\it Math. Ann. 290}, 671-680 (1991).

\bibitem {K}  J. Kazdan and F. Warner,  Curvature functions for open 2-manifolds, {\it Ann. Math. 99},  203-219 (1974).

\bibitem {McO1} R. McOwen, On the equation $\Delta u+Ke^{2u}=f$ and prescribed negative curvature in $\R^2$,
{\it J. Math. Anal. Appl. 103}, 365-370 (1984).

\bibitem {McO} R. McOwen, Conformal metrics in $\R^2$ with prescribed Gaussian curvature and positive
total curvature, {\it Indiana Univ. Math. J. 34}, 97-104  (1985).

\bibitem {N} W.-M. Ni,   On the elliptic equation $\Delta u+K(x)e^{2u}=0$ and conformal metric with
prescribed Gaussian curvatures, {\it Invent. Math. 66}, 343-352 (1982).

 \bibitem {S} D. Sattinger,   Conformal metrics in $\R^2$ with prescribed curvature,
 {\it Indiana Univ.  Math. J. 22}, 1-4 (1972).
\end{thebibliography}
\end{document}